\RequirePackage{fix-cm}
\documentclass[smallextended,natbib]{svjour3}       
\smartqed
\usepackage{amsmath}
\usepackage{amsfonts}
\usepackage{txfonts}
\usepackage{bm}
\usepackage{braket}
\usepackage{caption}
\usepackage{mathrsfs}
\usepackage{dcolumn}
\usepackage{graphicx}
\usepackage{epstopdf}
\usepackage[colorlinks=true,linkcolor=black,anchorcolor=black,citecolor=black,filecolor=black,menucolor=black,runcolor=black,urlcolor=blue]{hyperref}
\usepackage{dcolumn}
\usepackage[noend]{algpseudocode}
\allowdisplaybreaks

\DeclareMathOperator{\res}{Res}
\newcommand\Ord{\Or}

\newcommand\vecp[1]{\vec{#1}}

\newcommand\Lo{\hat{\boldsymbol{\mathsf L}}}
\newcommand\Po{\hat{\boldsymbol{\mathsf  P}}}
\newcommand\So{\hat{\boldsymbol{\mathsf  S}}}
\newcommand\R{\hat{\boldsymbol{\mathcal R}}}
\newcommand\Sr{\hat{\boldsymbol{\mathsf R}}}
\newcommand\D{\hat{\boldsymbol{\mathsf  D}}}
\newcommand\Rz{\hat{\boldsymbol{\mathcal Q}}}

\newcommand\Ha{{\partial H}/{\partial\ale}}
\newcommand\rmi{i}
\newcommand\rme{e}
\newcommand\rmd{d}
\newcommand{\Or}{\mathord{\mathrm{O}}} 
\newcommand{\eref}[1]{(\ref{#1})} 
\newcommand\exppe[2]{{\rm exp}_{+}(\phantom{-} \int_0^{#2}{#1}\,\rmd s)}
\newcommand\expme[2]{{\rm exp}_{-}(-\int_0^{#2}{#1}\,\rmd s)}
\newcommand\Us{\hat{\boldsymbol{\mathsf  U}}\vphantom{\boldsymbol{\mathsf  U}}}
\newcommand\Ui{\Us^{-1}}
\newcommand\expp[1]{\exppe{\Lo_{#1}}\ale}
\newcommand\expm[1]{\expme{\Lo_{#1}}\ale}
\newcommand\ale{\varepsilon}
\newcommand{\G}{W}
\date{\today}
\journalname{CeMDA}

\begin{document}
\sloppy
\title{Generalisation of the explicit expression for the Deprit generator to  Hamiltonians nonlinearly dependent on small parameter.}
\titlerunning{Generalisation of the explicit expression for the Deprit generator}
\author{Nikolaev~A.\,S.}
\institute{A. Nikolaev\at Institute of Computing for Physics and Technology, Protvino, Moscow reg., Russia, and
\at RDTeX LTD, Moscow, Russia\\
\url{http://andreysnikolaev.wordpress.com}\\
\email{Andrey.Nikolaev@rdtex.ru}
}
\maketitle

\begin{abstract}
This work explores a  structure of 
the Deprit perturbation series 
and its  connection to a Kato resolvent expansion. 
It extends the formalism previously developed for the Hamiltonians linearly dependent on  perturbation parameter to a nonlinear case.

We  construct a canonical intertwining of perturbed and unperturbed averaging operators. 
This  leads to an  explicit   expression  for the  generator of the Lie-Deprit transform in any perturbation order.
Using this  expression, we  discuss   a regular pattern in the series, 
non-uniqueness of the  generator and  normalised  Hamiltonian,
 and  the uniqueness of the Gustavson integrals.

Comparison  of the corresponding  computational algorithm with  
classical perturbation methods  demonstrates 
its competitiveness for  Hamiltonians with a limited number of perturbation terms.

\end{abstract}

\section{Introduction}

  In a previous article  \citep{nikolaev4}, we  demonstrated
   a connection between the canonical perturbation series in classical mechanics and a Kato  expansion for the resolvent of the Liouville operator.
However, only the perturbations having the particular  form $H_0+\ale H_1$  were considered. 

 In this paper we extend the formalism to general perturbed Hamiltonians  
 represented by the power  series 
$ H=H_0 + \sum_{0}^\infty \ale^k H_k$.

Our method was inspired by the remarkable analogies between mathematical formalisms of perturbation expansions in classical and quantum mechanics.
It uses    {Kato} series for  
 the Laurent coefficients of the 
  resolvent of the Liouville operator  \citep{Kato}. 
The approach reveals a  regular structure in the  perturbation series and 
leads to  an  explicit expression  for the  generator of a normalising transform in any perturbation order:
\[\G = \So_{H} \frac{\partial H}{\partial \ale}.\]
Here, $\So_H$  is the  partial pseudo-inverse of the  perturbed Liouville operator.
    The canonical intertwining of perturbed and unperturbed averaging operators  allows for a description of  
  non-uniqueness of the generator and normalised  Hamiltonian. 
Surprisingly,   Gustavson integrals are unique. 

After a description of the explicit computational algorithm, we will compare 
its results and efficiency 
 with those of the major canonical perturbation approaches by
  \citet{Gustavson},  Hori-\citet{Mersman}, \citet{Deprit}, \citet{Henrard72},
and \citet{DragtFinn}.
Our  experiments with high-order computations for
 simple models demonstrate
the competitiveness of the  explicit  algorithm for  Hamiltonians with a limited number of perturbation terms.

Note that all calculations  are  formal, in the sense that neither a discussion of power
series convergence nor conditions for the existence of constructed operators are present.

The downloadable 
supplementary data files \citep{Demo} contain demonstrations  and large formulae.
The demonstrations use the freeware computer algebra system FORM \citep{Form}.

\section{Basic perturbation operators} \label{sec:basic}
We consider a  Hamiltonian $H$
that differs from an exactly solvable one by a  power series perturbation:
\[ 
H=H_0+V(\ale),\qquad V = \sum_{k=1}^\infty \ale^k H_k.
\]
Here $H_0$ and $H_k$ are the functions of canonical variables
${\bf x} = (\vecp p, \vecp q)$ on 
$\mathbb{R}^{2{\textbf d}}$ phase space:
\[x_i=q_i,\qquad x_{i+{\rm\textbf d}}=p_i,\qquad i=1,\,\ldots\,
{\rm\textbf d},
\]
equipped with the canonical structure  (Poisson brackets):
\[
[F({\bf x}),G({\bf x})] =  \sum_{i=1}^{\rm\textbf d}
\frac{\partial F}{\partial q_i}\frac{\partial G}{\partial p_i} -
\frac{\partial F}{\partial p_i}\frac{\partial G}{\partial q_i}. 
\]

We will discuss  
autonomous (time-independent) Hamiltonians having compact energy surfaces $H(\vecp p,\vecp q)=E$.
We will also  assume that all  functions  are analytic. 

The goal of the canonical perturbation theory is to transform  the perturbed Hamiltonian into a simpler operator 
using a near-identity canonical (i.e. preserving   the Poisson brackets)  transform:
\[ {\bf x} = \Us(\ale)\,  \widetilde{\bf x},\qquad\Us \left([F,G]\right)=[\Us F,\Us G].\]
The families of 
such transforms may be  handled conveniently using an operator formalism \citep{Deprit,Giorgilli:arXiv1211.5674}. 
We will  outline  some of the basics of it.

\paragraph{Liouvillian:}
For any analytic $\G({\bf x})$, one may introduce  a  Liouville operator  acting in the  space of analytic functions:
\[
\Lo_{\G} = [ \,.\,,\G], \qquad \Lo_{\G}  F = [F, \G]. 
\]
Due to the  Jacobi identity,  the Liouvillian is the \emph{derivation} of a Poisson bracket: 
\begin{equation}
\Lo_{\G} [ F,  G]= [\Lo_{\G}  F, G ]+[  F,\Lo_{\G}  G]. \label{Jacobi}
\end{equation}

\paragraph{Averaging operator:}
We assume that the unperturbed Hamiltonian system with the compact energy surface $H_0({\bf{x}})=E$ is completely integrable in a Liouville sense and performs   quasi-periodic motion on the invariant tori  \citep{Arnold}.  
The time average
\begin{equation}
\langle F\rangle =
\lim_{T\to\infty}\frac{1}{T} \int_0^{T} {\left. F(\vecp p(t),\vecp q(t))\right|}_{p(0)=p\atop q(0)=q} \,\rmd t = 
\lim_{T\to\infty}\frac{1}{T} \int_0^{T}
\rme^{t\Lo_{H_0}} F({\bf{x}}) \rmd t\label{Cesaro}
\end{equation}
exists in an \emph{action-angle} representation for any analytic function $F(\bf{x})$.
This average  is a function of the initial point ${\bf x}$.
Because it is written in an invariant form, it also exists in  other canonical variables.

The time averaging operation extracts from $F(\bf{x})$ its  \emph{secular} non-oscillating part, which
remains constant under the time evolution.
In functional analysis, it is known as  {\em Ces\`aro $(C,1)$}   average. 

However, we  will use  a  stronger  {\em Abel} average which is  common in quantum physics:
\[
\langle F\rangle^{(A)} =
\lim_{\lambda\to +0}\lambda\int_0^{+\infty}\rme^{-\lambda t} \rme^{t\Lo_{H_0}} F({\bf x})\, \rmd t.
\]
Whenever the Ces\`aro  average \eref{Cesaro} exists, the  Abel averaging produces the same results \citep{HP}. 
The Abel average greatly simplifies the formulae  and leads to a natural  connection to resolvent formalism. 
Strictly speaking, we should discuss corresponding Tauberian theorems, 
but we have limited our focus   to  formal expressions only.  

The corresponding  averaging operator:
\begin{equation}
\Po_{H_0} = \lim_{\lambda\to +0}\lambda\int_0^{+\infty}\rmd t \,\rme^{-\lambda t} \rme^{t\,\Lo_{H_0}},\label{Pdef}
\end{equation}
 is a projector $\Po_{H_0}^2 = \Po_{H_0}$. It projects $F({\bf x})$ onto the  space of functions commuting with $H_0$. 
This is the kernel of $\Lo_{H_0}$ and the algebra of unperturbed integrals of motion.

The notation $\Po_{H_0}$ was borrowed from  quantum mechanics. Other common notations for this operator are  ${\pmb{\langle\ \rangle}}$, $\bf\overline{F}$ and $M_t$.

Classical perturbation theory provides several recipes  to compute the averaging operator.
In the  action-angle  coordinates, the unperturbed Hamiltonian 
is the function of   ${\textbf d}$ action variables $\vecp J$, and the perturbation is  $2\pi$ periodic in the phases $\vecp\phi$:
\begin{displaymath}
H = H_0(\vecp J) + V(\vecp J,\vecp\phi,\ale) 
= H_0(\vecp J) +  \sum_{\vecp k} \tilde V(\vecp J,\vecp k,\ale)\,  \rme^{\rmi(\vecp k,\vecp\phi)}.
\end{displaymath}
The motion of the unperturbed system is quasi-periodic:
$\vecp J = \textrm{const}$, $\vecp\phi(t) = \vecp\omega t+\vecp\phi_0$ and $\vecp\omega={\partial H_0}/{\partial\vecp J}$. 
Since the Fourier components in the action-angle representation are eigenfunctions of 
$\Lo_{H_0}$, 
 the result of the averaging operation can be  written in the well known form:
\begin{equation}
\Po_{H_0} F = \lim_{\lambda\to +0} \lambda \int_0^{+\infty} 
\rme^{-\lambda t}\sum_{\vecp k} \tilde F(\vecp J,\vecp k) \,
\rme^{\rmi(\vecp k,\vecp\phi_0)+ \rmi(\vecp\omega,\vecp k)t}\rmd t 
= \sum_{(\vecp\omega,\vecp k)=0} \tilde F(\vecp J,\vecp k) \,
\rme^{\rmi(\vecp k,\vecp\phi_0)}.\label{AAdefs}
\end{equation}

An another example is the  Birkhoff-Gustavson representation  for power series Hamiltonians \citep{Gustavson}. 
In the simplest case, a quadratic 
unperturbed Hamiltonian is diagonalisable into
$H_0 = \sum {\omega_k}\, (p_k^2+q_k^2)/{2}$.

 After the canonical transform to complex variables: 
\[q_k =  1/{\sqrt{2}}\,(\zeta_k + \rmi\, \eta_k), \qquad p_k = {\rmi}/{\sqrt{2}}\, (\zeta_k - \rmi\, \eta_k),
\qquad k=1,\ldots,{\textbf d},\]
the Hamiltonian becomes:
\[H=\sum_{k=1}^{\textbf d} \rmi\,\omega_k \zeta_k \eta_k +  \sum_{|\vecp m|+|\vecp n|\ge 3} \tilde V(\vecp m,\vecp n,\ale) \prod_{k=1}^{\textbf d} \zeta_k^{m_k} \eta_k^{n_k}.
\]
The monomials $\zeta^{\vecp m}\eta^{\vecp n}=\prod \zeta_k^{m_k} \eta_k^{n_k}$  are  eigenfunctions of the unperturbed Liouvillian:
\[\Lo_{H_0} \zeta^{\vecp m}\eta^{\vecp n} = \rmi\,(\vecp\omega,\vecp m - \vecp n)\, \zeta^{\vecp m}\eta^{\vecp n}.\]
 Therefore, for  any series $F({\vecp p,\vecp q})= \sum \tilde F(\vecp m,\vecp n)\,  \zeta^{\vecp m}\eta^{\vecp n}$,
 its formal average can be computed as:
\[
\Po_{H_0} F = \sum_{(\vecp\omega,\vecp m - \vecp n) =0} \tilde F(\vecp m,\vecp n)\, \zeta^{\vecp m}\eta^{\vecp n}.
\]

\paragraph{Integrating operator:}
The complementary projector $1-\Po_{H_0}$ extracts the time-oscillating part from $F({\bf x})$. 
It projects on the non-secular space of oscillating functions where the inverse of $\Lo_{H_0}$ exists
(we will consider  only
the semi-simple $\Lo_{H_0}$).

This inverse is called the integrating operator $\So_{H_0}$;  it is also called  the solution of the homological equation, tilde 
operation, zero-mean antiderivative,  
Friedrichs ${\bf \widehat\Gamma}$ operation, ${\bf{1}/{k}}$ operator and division 
operation. Its invariant definition was given by   \citet{Primas}: 
\begin{equation}
\So_{H_0} = - \lim_{\lambda\to{+0}}\int_0^\infty \rmd t \,\rme^{-\lambda t} \rme^{t\,\Lo_{H_0}} \left(1-\Po_{H_0}\right).\label{Sdef}
\end{equation}
In  the  action-angle representation, the result of the integrating operator is as follows:
\begin{equation}
\So_{H_0} F = \sum_{(\vecp\omega,\vecp k)\ne 0} \frac{1}{\rmi(\vecp\omega,\vecp k)}
\tilde F(\vecp J,\vecp k) \, \rme^{\rmi(\vecp k,\vecp\phi_0)}.\label{Sdef1}
\end{equation}
According to KAM theory \citep{Kolmogorov53}, for a general nondegenerate multidimensional
system, the operator \eref{Sdef1} is analytic for almost all frequencies except a set of Lebesgue measure zero. This
suffices for our formal constructions.

The corresponding expression in the Birkhoff-Gustavson representation is: 
\[
\So_{H_0} F = \sum_{(\vecp\omega,\vecp m - \vecp n)\ne 0} \frac{1}{\rmi\,(\vecp\omega,\vecp m - \vecp n)}
\tilde F(\vecp m,\vecp n)\, \zeta^{\vecp m}\eta^{\vecp n}.\nonumber
\]

The operators $\Lo_{H_k}$, $\Po_{H_0}$ and  $\So_{H_0}$ are  the building blocks of the canonical perturbation series presented here. It is easy to check that:
\begin{align*}
\Lo_{H_0}\So_{H_0} &= \So_{H_0}\Lo_{H_0}=1-\Po_{H_0},\\
\So_{H_0}\Po_{H_0} &= \Po_{H_0}\So_{H_0}\equiv 0.
\end{align*}
Further {canonical identities} have been discussed in  \citet{nikolaev4}. 

\section{Canonical perturbation theory} \label{sec:program}

It is  known  \citep{Arnold} that any near-identity family of canonical transforms 
\[{\bf x} = \Us_\G(\ale)\, \widetilde {\bf x},
\qquad\widetilde H=\Us_\G H,
\]
is a Hamiltonian flow  in time $\ale$ \ with some generator $\G({\bf x},\ale)$. The corresponding  operators obey the equations:
\begin{equation}
\frac{\partial\Us_\G}{\partial\ale} =  \Us_\G \Lo_\G,\qquad \frac{\partial\Ui_\G}{\partial\ale} = - \Lo_\G\Ui_\G .
\label{depriteq}
\end{equation}
In classical mechanics, such $\ale$-dependent parameterisations  
are called  Lie-Deprit transforms \citep{Deprit}.
In quantum physics,  corresponding objects  are known as  ordered exponentials  \citep{Dyson49}. 

When comparing  the quantum analogues \citep{Suzuki85,Nikolaev5}, we should mention that 
the classical Lie-Deprit transform $\Us_\G$ corresponds to the
 negative 
 ordered quantum exponential $\expm{(-\G)}$,  the inverse transform $\Ui_\G$ corresponds to 
 the positive ordered exponential $\expp{(-\G)}$ and the generator has the opposite sign.

The perturbation theory constructs the transforms and the generator in the form of  a power series:
\[\Us_\G=\sum_{n=0}^\infty \ale^n\, \Us_n,\qquad \Ui_\G=\sum_{n=0}^\infty \ale^n\, \Ui_n,\qquad \G({\bf x},\ale) = \sum_{n=0}^\infty \ale^n \G_n({\bf x}).\]
To demonstrate the regular structure of the perturbative series, we do not 
include  factorial denominators in the above expressions, which is traditionally done.
To avoid any misunderstandings, we should note that $\Ui_n$ means the operator coefficient of $\ale^n$ in the series for $\Ui_\G$.
 This is not the inverse of the $\Us_n$ coefficient.

The substitution of the series for $\G$, $\Us$ and $\Ui$ into  the equations \eref{depriteq}
results in the following recursive relations for the coefficients: 
\begin{equation}
\Us_n =  \frac{1}{n} \sum_{k=0}^{n-1}  \Us_k \,\Lo_{\G_{n-k-1}}, \qquad
\Ui_n = - \frac{1}{n} \sum_{k=0}^{n-1}  \Lo_{\G_{n-k-1}} \Ui_k.\label{Caryrel}
\end{equation}
Iterating  these expressions,  \citet{Cary} obtained the   non-recursive formulae: 
\begin{align*}
\Us_n^{\phantom{-1}}&= \sum_{\substack{(m_1,\ldots,m_r)\\n>m_1>m_2>\cdots>m_r}} \;\;\;\;\;\;
\frac{\Lo_{\G_{m_r-1}}}{m_r} \cdots\frac{\Lo_{\G_{m_1-m_2-1}}}{m_1}\frac{\Lo_{\G_{n-m_1-1}}}{n}\ ,\\
\Ui_n&= \sum_{\substack{(m_1,\ldots,m_r)\\n>m_1>m_2>\cdots>m_r}} 
\!\!\!\!\!\!(-1)^{r+1}
\frac{\Lo_{\G_{n-m_1-1}}}{n} \frac{\Lo_{\G_{m_1-m_2-1}}}{m_1}\cdots\frac{\Lo_{\G_{m_r-1}}}{m_r}\ .
\end{align*}
The sum runs  over all sets of 
integers $(m_1,\ldots,m_r)$, satisfying $n>m_1>\cdots>m_r>0$.
In this paper, we will derive a similar formula for the normalising  generator.

The canonical perturbation theory  constructs $\G= \sum_0^\infty \ale^n \G_n({\bf x})$  to simplify the transformed Hamiltonian as much as possible.
In the first order:
\begin{align*}\widetilde H&=  \Us_\G (H_0 + \sum_{k=0}^\infty \ale^k H_k) =H_0+\ale (H_1+\Lo_{\G_0}H_0) +\Ord(\ale^2)\\
&=H_0+\ale\left(H_1-\Lo_{H_0} \G_0 \right)+\Ord(\ale^2).
\end{align*}

If one choose $\G_0=\So_{H_0}H_1$ in the above expression, then all the terms of  order  $\ale$  in the  transformed Hamiltonian  
become secular (beginning with $\Po_{H_0}$ operator): 
\begin{equation*}
\widetilde H = H_0+\ale \,\Po_{H_0} H_1 +\Ord(\ale^2).
\end{equation*}

Typically, a perturbation method constructs the normalising generator and the transformed Hamiltonian order by order.
It  consequentially chooses $\G_1$ 
to eliminate non-secular terms up to $\ale^2$, then  $\G_2$  to eliminate nonsecular terms up to $\ale^3$, and so on.
The detailed description of such procedures can be found in classical books on the perturbation theory \citep{Naifeh,Giacaglia,FM}.

Obviously, the  generator $\G$ is not unique. An arbitrary secular function $\Po_{H_0} F$ may be added to any order. Further orders  
will depend on this term. The rule for choosing the secular part of $\G$ is the   {\em uniqueness condition}. 
It is common to construct a completely non-secular  generator   
$\Po_{H_0} \G_{NS} = 0$.
However, we will see that other conditions may be useful.

If the process converges, the  canonical transform with the generator $\G=\sum_0^\infty \ale^n \G_n$  
  makes  the  Hamiltonian completely secular
\[\Po_{H_0} \widetilde H = \widetilde H,\]
 or  equivalently, $\Ui_\G \, \Po_{H_0} \Us_\G H =  H$. 

Conversely, if we  construct  an analytic  continuation of 
the averaging operator $\Po_{H_0}$
 to the perturbed case
\[\Po_H(\ale) H(\ale) = H,\qquad \Po_H(0)=\Po_{H_0},\]
and a transform,
 such that the perturbed and unperturbed projectors are 
 canonically connected
\begin{equation}
 \Po_H = \Ui_\G \, \Po_{H_0} \Us_\G,\label{Symeq}
 \end{equation}
 then, due to the intertwining properties of  this transform
\begin{equation}
 \Us_\G\, \Po_H = \Po_{H_0}\, \Us_\G,\qquad \Po_H \Ui_\G= \Ui_\G \Po_{H_0},\label{Intereq}
 \end{equation}
  a new Hamiltonian will become  secular: 
\begin{equation*}
\Po_{H_0} \widetilde H = \Po_{H_0}  \Us_\G H = 
 \Us_\G \Po_H H 
=  \Us_\G H = \widetilde H.
\end{equation*}
Therefore, this construction will explicitly realise the goal of the canonical perturbation theory. 
In the next sections, we will  develop such a continuation and  intertwining canonical transform
using the methods of functional analysis.

\section{The Liouvillian resolvent}
Consider the resolvent of the Liouville operator:
\begin{equation}
\R_{H}(z) =\left({{\Lo_{H}-z}}\right)^{-1}.\label{Res}
\end{equation}
This  operator-valued function of the complex variable $z$ is the Laplace transform of  the evolution operator   of the Hamiltonian system:
\[
\R_{H}(z) = - \int_0^{+\infty}\rmd t\, \rme^{-z t} \rme^{t\Lo_{H}}.
\]

Resolvent singularities are the eigenvalues of   $\Lo_{H}$. For an integrable Hamiltonian system with compact energy surfaces, these eigenvalues belong to an imaginary axis. 
Typically, the  spectrum of a Liouville operator is anywhere dense \citep{spohn75}. 

Let us begin with a simple case of an isolated point spectrum. We  consider a one-dimensional  system and 
 restrict the  domain of the resolvent operator to  analytic functions 
with the argument  on the compact energy surface $H({\bf x})=E$.
Under such conditions, the system is non-relaxing and oscillates with the single frequency $\omega(E)$. The resolvent singularities are located  at the points 0, $\pm \rmi\omega(E)$, $\pm 2\rmi\omega(E)$, \ldots .

The existence of  the operator $\Po_{H_0}$ 
defined by \eref{Pdef} 
means that the unperturbed resolvent has a simple pole at $0$. 
The averaging operator is the residue of the resolvent at that pole:
\[
\Po_{H_0} \equiv - \res_{z=0} \R_{H_0},
\]
while the integrating operator $\So_{H_0}$ is its holomorphic part:
\[
\So_{H_0}=\lim_{z\to0}{\R_{H_0}(z)\left(1-\Po_{H_0}\right)}.
\] 
Therefore, the Liouvillian resolvent combines both of the basic perturbation operators  \citep{nikolaev4}.
This allows us to apply the  powerful  formalism of complex analysis to the perturbation theory.  

It is well known \citep{Kato} that, due to the Hilbert resolvent identity: 
\begin{equation}
\R_{H}(z_1)-\R_{H}(z_2)=(z_1-z_2)\,\R_{H}(z_1)\R_{H}(z_2),\label{HilbertId}
\end{equation}
the Laurent series for the unperturbed resolvent  has the  form of a geometric progression:
\begin{equation}
\R_{H_0(z)}=-{1\over z}\Po_{H_0} + \sum_{n=0}^\infty z^n \So_{H_0}^{n+1}=\sum_{n=0}^{+\infty} \Sr_{H_0}^{(n)}z^{n-1}.\label{UResLoran}
\end{equation}
Here, we  denoted $\Sr_{H_0}^{(0)}=-\Po_{H_0}$ and $\Sr_{H_0}^{(n)}= \So_{H_0}^n$.

The perturbed resolvent  may be more singular. The Laurent series for a general resolvent with an isolated singularity at the origin has
the following form \citep{Kato}:
\[
\R_H(z)=-{1\over z}\Po_H + \sum_{n=0}^\infty z^n \So_H^{n+1}
-\sum_{n=2}^\infty z^{-n} \D_H^{n-1}=\sum_{n=-\infty}^{+\infty} \Sr_{H}^{(n)}z^{n-1},\label{ResLoran}
\]
where $\D_H$ is the eigen-nilpotent operator, which does not have an unperturbed analogue ($\D_{H_0} \equiv 0$). 

\subsection{Kato series}
In an our case of the isolated point spectrum, the perturbed resolvent can be expanded into  the Neumann series.
Since
\[
 \Lo_{H_0 + V}-z=\left(1+\Lo_{V}\R_{H_0}\right)\left(\Lo_{H_0}-z\right),
\]
 then for $V = \sum_{k=1}^\infty \ale^k H_k$, the following expansion holds:
\begin{equation}
\R_{H}(z)=\R_{H_0}{\left(1 +\Lo_{V}\R_{H_0}\right)}^{-1}
=\R_{H_0}(z)\sum_{n=0}^{\infty}{\left(-\Lo_{V}\R_{H_0}(z)\right)}^n
=\R_{H_0}(z)+\sum_{n=1}^{\infty}\ale^n\, \R_{H}^{(n)}(z).\label{Neumann}
\end{equation}
The coefficient operators 
are given by the  expression:
\[\R_{H}^{(n)}(z)=\sum_{m=1}^n (-1)^m \sum_{\substack{
k_1+\cdots + k_m={ n}\\ k_j> 0}} 
 \R_{H_0}(z)\,\Lo_{H_{k_m}}\R_{H_0}(z)\,\ldots
\Lo_{H_{k_1}}\R_{H_0}(z),
\]
where the sum being taken for all combinations of $1\le m\le n$ positive integers $k_1,\ldots k_m$, such that 
$k_1+\cdots + k_m={ n}$.

The integration  
around a sufficiently small contour $|z|=a$ 
encompassing the singularity at the origin results in the Kato   series 
 for the \emph{perturbed averaging operator} \citep{Kato}:
\begin{align*}
\Po_{H} &= -{1\over{2\pi \rmi}}\oint_{|z|=a} \R_{H}(z)\rmd z\\
&=\Po_{H_0}-{1\over{2\pi \rmi}}\sum_{n=1}^\infty\ale^n\oint_{|z|=a}\sum_{m=1}^n (-1)^m \sum_{\substack{
k_1+\cdots + k_m={ n}\\ k_j\ge 1}} 
 \R_{H_0}(z)\,\Lo_{H_{k_m}}\R_{H_0}(z)\,\ldots
\Lo_{H_{k_1}}\R_{H_0}(z) \, \rmd z\\
&= \Po_{H_0}-{1\over{2\pi \rmi}}\sum_{n=1}^\infty\ale^n\oint_{|z|=a}\sum_{m=1}^n (-1)^m \sum_{\substack{k_1+\cdots + k_m={ n}\\ k_j\ge 1}} 
\left(\sum_{p_{m+1}=0}^\infty \Sr_{H_0}^{(p_{m+1})}z^{p_{m+1}-1}\right)\,\Lo_{H_{k_m}}
\\
&\qquad\qquad\qquad
\,\left(\sum_{p_{m}=0}^\infty \Sr_{H_0}^{(p_{m})}z^{p_{m}-1}\right)\ \ldots
\ \Lo_{H_{k_1}}\left(\sum_{p_{1}=0}^\infty \Sr_{H_0}^{(p_{1})}z^{p_{1}-1}\right) \, \rmd z.
\end{align*}
 The \emph{perturbed integrating operator} and the \emph{perturbed eigen-nilpotent operator} are:
\begin{equation*}
\So_H = {1\over{2\pi \rmi}}\oint_{|z|=a} z^{-1}\R_H(z)\rmd z,\qquad \D_H = -{1\over{2\pi \rmi}}\oint_{|z|=a} z\,\R_H(z)\rmd z .
\end{equation*}
Only the coefficients of $z^{-1}$ in the above  expressions will contribute to the  result:
\begin{align}
\Po_H&=\Po_{H_0}+\sum_{n=1}^\infty\ale^n\left(\sum_{m=0}^n {(-1)}^{m+1}\sum_{\substack{
 p_1+\cdots + p_{m+1}={m}\\k_1+\cdots + k_m={ n\hphantom{+1}}\\ k_j\ge 1, p_j \ge0}}
\Sr_{H_0}^{(p_{m+1})}\Lo_{H_{k_m}}\Sr_{H_0}^{(p_m)}\,\ldots\,\Sr_{H_0}^{(p_2)}
\Lo_{H_{k_1}}\Sr_{H_0}^{(p_1)}\right),\nonumber\\
\So_H&=\So_{H_0}+\sum_{n=1}^\infty\ale^n\left(\sum_{m=0}^n {(-1)}^{m\hphantom{+1}}\sum_{\substack{
 p_1+\cdots + p_{m+1}={m+1}\\k_1+\cdots + k_m={ n}\\ k_j\ge 1, p_j \ge0}}
\!\!\!\!\Sr_{H_0}^{(p_{m+1})}\Lo_{H_{k_m}}\Sr_{H_0}^{(p_m)}\,\ldots\,\Sr_{H_0}^{(p_2)}
\Lo_{H_{k_1}}\Sr_{H_0}^{(p_1)}\right),\nonumber\\
\D_H&=\phantom{\So_{H_0}+{}}\sum_{n=1}^\infty\ale^n\left(\sum_{m=0}^n {(-1)}^{m+1}\sum_{\substack{
 p_1+\cdots + p_{m+1}={m-1}\\k_1+\cdots + k_m={n}\\ k_j\ge 1, p_j \ge0}}
\!\!\!\!\Sr_{H_0}^{(p_{m+1})}\Lo_{H_{k_m}}\Sr_{H_0}^{(p_m)}\,\ldots\,\Sr_{H_0}^{(p_2)}
\Lo_{H_{k_1}}\Sr_{H_0}^{(p_1)}\right).\label{shph}
\end{align}

Because the operator $\Sr_{ H_0}^{(p_j)}$ consists of $p_j$ operators  $\So_{H_0}$,
the summation in the coefficient of $\ale^n$ in the above expressions is being taken by all possible 
subdivisions of $n$ in $m$ parts by perturbations $\Lo_{H_{k_m}}$ and all placements  of  $m$ (or $m+1$, or  $m-1$)  
integrating operators
$\So_{H_0}$ in these parts. 
For example, the first two orders of the perturbed integrating operator are:
\begin{align*}
\So_H\!=&\So +\!\ale\left( \So^2\Lo_{H_1}\Po\! +\! \Po\Lo_{H_1}\So^2\!\!-\!\So\Lo_{H_1}\So\right) 
+\!\ale^2\left(\So^2\Lo_{H_2}\Po\!-\!\So\Lo_{H_2}\So\!+\!\Po\Lo_{H_2}\So^2\!\!+\!\So\Lo_{H_1}\So\Lo_{H_1}\So\right.\\
&{}-\So^2\Lo_{H_1}\So\Lo_{H_1}\Po- \So^2\Lo_{H_1}\Po\Lo_{H_1}\So-\Po\Lo_{H_1}\So\Lo_{H_1}\So^2
- \Po\Lo_{H_1}\So^2\Lo_{H_1}\So-\So\Lo_{H_1}\So^2\Lo_{H_1}\Po  \\
&{}- \So\Lo_{H_1}\Po\Lo_{H_1}\So^2+\Po\Lo_{H_1}\Po\Lo_{H_1}\So^3
+\Po\Lo_{H_1}\So^3\Lo_{H_1}\Po +\left.\So^3\Lo_{H_1}\Po\Lo_{H_1}\Po\right) + \Ord(\ale^3).
\end{align*}
Here we  omitted the subscript $H_0$ for  the {unperturbed} $\Po_{H_0}$ and $\So_{H_0}$ operators.

The properties of the unperturbed operators can be  extended to their analytic continuations  
as follows:
\begin{align}
&\Po_{H} {H} = { H}, \qquad\So_{H}\Lo_{H} =\Lo_{H}\So_{H}= 1-\Po_{H},\qquad\Lo_{ H} \Po_{H} = \Po_{H} \Lo_{H} = \D_{H}, \nonumber\\
&\So_{H}\Po_{H}=0,
\qquad\! \Po_{H} \D_{H}=\D_{H} \Po_{H}=\D_{H}. 
\label{pertprop}
\end{align}
For the details, refer to  \citet{nikolaev4}.

To avoid any misunderstandings, it should be noted that $\Po_{H} F$ will not be  commuting with  the perturbed Hamiltonian. 
This is because, in general, $\Lo_{H} \Po_{H} =  \D_{H}\neq 0$. 
Actually, the  perturbed projector,  $\Po_{H}$, projects onto some analytic continuation
of the algebra of the integrals of the unperturbed Hamiltonian. 
This may not coincide with the algebra of the integrals of the perturbed system because of a destruction of  symmetries.

\subsection{Canonical properties of the  Liouvillian resolvent}

Since  $\Lo$ is a derivative, there  exists an  { integration by parts} formula for its resolvent. 
It results from the  application  of the identical operator
\[
\R_{H}(z_1)\,\left(\Lo_{H} - z_1\right) \equiv  1,
\] 
 to the Poisson bracket $\Lo_{\,\R_{H}(z_2) \,F}\,\R_{ H}(z_3)$.
We suppose that complex $z_1$, $z_2$, and $z_3$ are outside of the spectrum  of $\Lo_{H}$ and that an arbitrary function $F({\bf{x}})$  is analytic.
 
An expansion of the Jacobi identity 
$\Lo_F\Lo_G - \Lo_G\Lo_F = \Lo_{\Lo_F G}$
 yields:
\begin{align*}
\Lo_{\R_{H}(z_2) \,F}\,\R_{ H}(z_3)
={}&\R_{H}(z_1)\,\Lo_{H} \,\Lo_{\R_{H}(z_2) \,F}\,\R_{ H}(z_3) - z_1  \R_{H}(z_1)\, \Lo_{\R_{H}(z_2) \,F}\,\R_{ H}(z_3) \\
={}&\R_{H}(z_1)\, \Lo_{\R_{H}(z_2) \,F}\Lo_{H}\,\R_{ H}(z_3) +  \R_{H}(z_1)\,\Lo_{\Lo_H {\R_{H}(z_2) \,F}}\,\R_{ H}(z_3)   \\
&\qquad- z_1  \R_{H}(z_1)\, \Lo_{\R_{H}(z_2) \,F}\,\R_{ H}(z_3)\\
={}&\R_{H}(z_1)\, \Lo_{\R_{H}(z_2) \,F}\left(1+z_3\R_{H}(z_3)\right) +  \R_{H}(z_1)\,\Lo_{\left(1+z_2\R_{H}(z_2)\right) \,F}\,\R_{ H}(z_3)   \\
&\qquad- z_1  \R_{H}(z_1)\, \Lo_{\R_{H}(z_2) \,F}\,\R_{ H}(z_3).
\end{align*} 
Therefore,  the following   identity holds true for the Liouvillian resolvent:
\begin{align}
 \R_{H}(z_1)\, \Lo_{F} \, \R_{H}(z_3)=&\Lo_{\R_{H}(z_2)\, F} \,\R_{H}(z_3)-\R_{H}(z_1) \,\Lo_{\R_{H}(z_2) \, F} \nonumber\\
&+ (z_1-z_2-z_3)\, \R_{H}(z_1) \,\Lo_{\R_{H}(z_2)\, F}\,  \R_{H}(z_3).\label{SimpID1}
\end{align}

Consider the derivative of the  resolvent  with respect to the perturbation
\[
\frac{\partial}{\partial\ale}\left(\R_{ H}(z)\right) = -\R_{ H}(z)\,\Lo_{\Ha}\,\R_{ H}(z).
\]
Substituting $z_1=z_3=z$ and $ F=\Ha$  in the  canonical  
identity \eref{SimpID1} 
 yields:
\vbox{\begin{align}
\R_{ H}(z)\,\Lo_{\,\Ha}\,\R_{ H}(z) =& \,\Lo_{\R_{ H}(z_2)\, \Ha}\,\R_{ H}(z) - \R_{ H}(z)\,\Lo_{\R_{ H}(z_2)\,\Ha} \nonumber\\
&- z_2\,  \R_{ H}(z)\, \Lo_{\R_{ H}(z_2)\, \Ha}\,\R_{ H}(z). \label{SimpID2}
\end{align}}
The coefficient of $z_2^0$ in the Laurent series for this  expression is the following:
\[
\frac{\partial}{\partial\ale}\R_{ H}(z) = 
\R_{ H}(z)\Lo_{\So_{ H}  \,\Ha} - \Lo_{\So_{ H}\,\Ha}\R_{ H}(z) -  \R_{ H}(z)\,\Lo_{\Po_{ H}\,\Ha}\,\R_{ H}(z).
\]
Proceeding in a similar way for the coefficients of $z_2^{-n}$ \hbox{$(n\ge1)$} in \eref{SimpID2}, 
we obtain: 
\begin{align*}
\R_{ H}(z)\Lo_{\Po_{ H} \,\Ha}&=\Lo_{\Po_{ H} \,\Ha}\,\R_{ H}(z) - \R_{ H}(z)\,\Lo_{\D_{ H} \,\Ha}\,\R_{ H}(z),\\
\R_{ H}(z)\Lo_{\D^n_{ H}\, \Ha}&=\Lo_{\D^n_{ H} \,\Ha}\,\R_{ H}(z) - \R_{ H}(z)\,\Lo_{\D^{n+1}_{ H} \Ha}\,\R_{ H}(z).
\end{align*}
This allows  us to rewrite the expression for  the  resolvent derivative as:
\begin{align*}
{\partial\over \partial\ale}\R_{ H}(z) =&
\,\R_{ H}(z)\,\Lo_{\So_{ H}\, \Ha} - \Lo_{\So_{ H}\, \Ha}\,\R_{ H}(z) \\
&- \Lo_{\Po_{ H}\, \Ha}\R_{ H}(z)^2 
+  \Lo_{\D_{ H}\, \Ha}\R_{ H}(z)^3
- \Lo_{\D^2_{ H}\, \Ha}\R_{ H}(z)^4 + \ldots
\end{align*}
Actually, this is a power series because $\D^n_{ H} = \Ord(\ale^n)$.

It follows from the Hilbert resolvent identity that:
\[
{\partial^n\over \partial z^n}\,\R_{ H}(z) = n!\ \R_{ H}^{n+1}(z).
\]
Finally, we obtain the following series:
\begin{align}
\frac{\partial}{\partial\ale}\R_{ H}(z) =&
\,\R_{ H}(z)\,\Lo_{\So_{ H} \,\Ha} - \Lo_{\So_{ H} \,\Ha}\R_{ H}(z) 
- \Lo_{\Po_{ H} \,\Ha}\frac{\partial\,\R_{ H}(z)}{\partial\,z} \nonumber\\
&{}+ \frac{1}{2} \Lo_{\D_{ H}\, \Ha}\frac{\partial^2\R_{ H}(z)}{\partial\,z^2}
- \frac{1}{6} \Lo_{\D_{ H}^2\, \Ha}\frac{\partial^3\R_{ H}(z)}{\partial\,z^3} 
+ \ldots
\end{align}
The derivative of the projector ${\partial \,\Po_{ H}}/{\partial\ale}$ can be obtained as the residue of this expression at $z=0$.
In the  case of an isolated point spectrum,  the Liouvillian resolvent  is a meromorphic function and, therefore the residue of any of its derivatives with respect to $z$  vanishes identically. 
As a result, the projector $\Po_{H}$ transforms
canonically under the perturbation:
\begin{equation}
\frac{\partial}{\partial\,\ale}\Po_{ H} = \Po_{ H}\,\Lo_{\So_{ H}\, \Ha} - \Lo_{\So_{ H} \,\Ha}\,\Po_{ H},\label{Projeq}
\end{equation}
and the projectors are connected by the canonical Lie-Deprit transform:
\[
\Po_{ H} =\Ui_{\G}\, \Po_{H_0}\,\Us_{\G},\qquad \G=\So_{ H}\,  \Ha.
\]
This  holds for all perturbation orders.

\section{An explicit expression for a  generator}

\subsection{Regular pattern in perturbation series}

As was discussed,  the canonical connection of the projectors means that   the Lie-Deprit transform  
with the generator
\begin{equation}
\G=\So_{ H}\, \frac{\partial H}{\partial \ale}\label{exp1}
\end{equation}
formally normalises the Hamiltonian in all orders in $\ale$.  
The expanded form of this expression demonstrates the regular pattern:
\begin{equation}
\G=\sum_{n=0}^\infty\ale^{n} \left(\sum_{m=0}^{n} {(-1)}^{m}\sum_{\substack{
 p_0+\cdots + p_{m}={\bf m+1}\\k_0+\cdots + k_m={n+1}\\ k_j\ge 1, p_l \ge0}}
\!\!\!\! k_0 \, \Sr_{H_0}^{(p_{m})}\Lo_{H_{k_m}}\Sr_{H_0}^{(p_{m-1})}\,\ldots\,\Sr_{H_0}^{(p_1)}
\Lo_{H_{k_1}}\Sr_{H_0}^{(p_0)}H_{k_0}\right) . 
\label{explicit}
\end{equation}
Here, the sum being taken for all the combinations of positive integers $k_0,\ldots k_m$ and nonnegative $p_0,\ldots,p_{m}$, $1\le m\le n$ such that  $\sum p_j={m+1}$, and $\sum k_j={n+1}$.

The first two perturbative orders for the  generator are:
\begin{align*}
\G=&\,\So H_1+\ale\left( \So^2\Lo_{H_1}\Po H_1 -\So\Lo_{H_1}\So H_1 + 2\,\So H_2    
\right)
+ \ale^2 \left(\So\Lo_{H_1}\So\Lo_{H_1}\So H_1\right.      \\
&{}- \So\Lo_{H_1}\So^2\Lo_{H_1}\Po H_1 -\So^2\Lo_{H_1}\So\Lo_{H_1}\Po H_1
             - \So^2\Lo_{H_1}\Po\Lo_{H_1}\So H_1-\Po\Lo_{H_1}\So\Lo_{H_1}\So^2 H_1 \\
&{} - \Po\Lo_{H_1}\So^2\Lo_{H_1}\So H_1 + \Po\Lo_{H_1}\So^3\Lo_{H_1}\Po H_1+ \Po\Lo_{H_1}\Po\Lo_{H_1}\So^3 H_1
+3\,\So H_3  -  2\,\So\Lo_{H_1}\So H_2 
  \\
&{}- \So\Lo_{H_2}\So H_1+ 2\,\So^2\Lo_{H_1}\Po H_2  
+\So^2\Lo_{H_2}\Po H_1+ 2\,\Po\Lo_{H_1} \So^2 H_2 
\left. + \Po\Lo_{H_2} \So^2 H_1\right) 
          + \Ord(\ale^3).
\end{align*}
The four orders for the normalised  Hamiltonian are:
\begin{align*}
\widetilde H&=H_0 
       + \ale\, \Po H_1 
       + \ale^2 \left( - \frac{1}{2}\Po\Lo_{H_1}\So H_1 + \Po H_2 \right)\\  
&{}       + \ale^3 \left( \frac{1}{3} \Po\Lo_{H_1}\So\Lo_{H_1}\So H_1
                 - \frac{1}{6} \Po\Lo_{H_1}\So^2\Lo_{H_1}\Po H_1
                 - \Po\Lo_{H_1}\So H_2    +\Po H_3
          \right) 
          \\
&{}+\ale^4 \left(\frac{1}{6}\Po\Lo_{H_1}\So\Lo_{H_1}\So^2\Lo_{H_1}\Po H_1 \right.
                 {}- \frac{1}{4}  \Po\Lo_{H_1}\So\Lo_{H_1}\So\Lo_{H_1}\So  H_1 
                 + \frac{1}{12} \Po\Lo_{H_1}\So^2\Lo_{H_1}\So\Lo_{H_1}\Po H_1\\
&{}               + \frac{1}{8} \Po\Lo_{H_1}\So^2\Lo_{H_1}\Po\Lo_{H_1}\So H_1  
               + \frac{1}{4} \Po\Lo_{H_1}\Po\Lo_{H_1}\So^2\Lo_{H_1}\So H_1
{}                 + \frac{1}{4} \Po\Lo_{H_1}\Po\Lo_{H_1}\So\Lo_{H_1}\So^2 H_1\\
&{}                 - \frac{1}{6} \Po\Lo_{H_1}\Po\Lo_{H_1}\So^3\Lo_{H_1}\Po H_1 
  {}               - \frac{1}{4} \Po\Lo_{H_1}\Po\Lo_{H_1}\Po\Lo_{H_1}\So^3 H_1 
  + \frac{1}{2}\Po\Lo_{H_1}\So\Lo_{H_1}\So H_2\\
&{}  + \frac{1}{4}\Po\Lo_{H_1}\So\Lo_{H_2}\So H_1  
- \frac{1}{4}\Po\Lo_{H_1}\So^2\Lo_{H_1} \Po H_2 
- \frac{1}{12}\Po\Lo_{H_1}\So^2\Lo_{H_2}\Po H_1  
- \frac{1}{2}\Po\Lo_{H_1}\Po\Lo_{H_1}\So^2 H_2 \\
&{}  
- \frac{1}{4}\Po\Lo_{H_1}\Po\Lo_{H_2}\So^2 H_1 
+ \frac{1}{4}\Po\Lo_{H_2}\So\Lo_{H_1}\So H_1
 - \frac{1}{6}\Po\Lo_{H_2}\So^2\Lo_{H_1}\Po H_1
-\frac{1}{2}\Po\Lo_{H_2}\So H_2\\
&{}\left.- \Po\Lo_{H_1} \So H_3
+ \Po H_4 
\right)+ \Ord(\ale^5),
\end{align*}
 Here we   omitted the indices $H_0$ for the unperturbed operators $\Po_{H_0}$ and $\So_{H_0}$
for compactness and we took into account the identities $\Po_{H_0}\Lo_{F}\So_{H_0}^2 F \equiv 0$ and $\Po_{H_0}\Lo_F\Po_{H_0} F\equiv 0$ \citep{nikolaev4}. 
These expressions  are the generalisation of the classic result by \citet{Burshtein}.

It is possible to extend   formula \eref{explicit} to multidimensional  systems 
in the form of an asymptotic series. 
Similarly to   previous  work \citep{nikolaev4}, it can be shown  that the expression for $\G$ truncated at order $\Ord(\ale^{N})$   normalises  Hamiltonian up to $\Ord(\ale^{N+1})$.
For  clarity, we  here  use formal analytic  expressions. However, for multidimensional systems,
  these expressions should be converted straightforwardly  into 
 truncated sums,  and all  equalities hold up to   $\Ord(\ale^{N+1})$.

\subsection{General form of the generator}

Let us  determine the general form of a generator of
  a transform connecting the perturbed and unperturbed projectors. 
It follows from  $(\ref{Symeq})$ that $\Po_H$ satisfies the operatorial differential equation:
\[
\frac{\partial}{\partial\,\ale}\Po_H = \Po_H \Lo_\G -\Lo_\G \Po_H  .
\]
Application of  this expression  to the Hamiltonian $H$  and the identities
\[\Po_{ H} { H} = { H},\qquad  
{\partial\over \partial\,\ale}\left(\Po_{ H} { H}\right) =
 \left({\partial\over \partial\,\ale}\Po_{ H} \right)  {H} + 
\Po_{H}  \frac{\partial}{\partial\,\ale}  {H},
\]
yield
\[
\left(1-\Po_H\right)\, \Lo_H \, \G =  \left(1-\Po_H\right)\,  \frac{\partial H}{\partial\,\ale}.
\]
To solve this equation, it is sufficient to apply the $\So_H$ operator. Therefore, the general form of the generator of the connecting transform is:
\begin{equation}
\G =  \So_H \frac{\partial H}{\partial\,\ale} + \Po_H F,\label{Final}
\end{equation}
where $F({\bf x},\ale)$  may be any analytic function, and $\Po_H$ and $\So_H$ are given by \eref{shph}.

This formula provides the non-recursive expression for
the Deprit generator of a normalising transform and defines it unambiguously. 
It generalises the corresponding formula in \citet{nikolaev4}.

The  choice of  $\Po_{H} F({\bf x},\ale)$ is the {\em uniqueness  condition}. It is natural to set $F\equiv 0$ or 
\begin{equation}
\Po_{ H} \G = 0.\label{uniq}
\end{equation}
This is not equal to the non-secular condition $\Po_{H_0} \G_{NS} = 0$, which is  traditionally used  in  
the canonical perturbation theory. 
Because $\Lo_{\Po_{H} F}\, \Po_{H} = \Po_{H} \Lo_{\Po_{H} F}$, the projector $\Po_{H}$  is itself insensitive to this choice.

We can conclude that the generators of normalising transforms 
can differ by a function belonging to  a continuation of  the algebra of the integrals of the unperturbed system.

Hamiltonians that were  normalised using the different uniqueness conditions $F_1({\bf x},\ale)$ and $F_2({\bf x},\ale)$
 are connected by a canonical transform
$ \Us_{21} = \Us_{\G_2} \,\Ui_{\G_1}$.
Let us find its generator:  
\[ \frac{\partial}{\partial\ale}\Us_{21} = \Us_{\G_2}
\ \left(\Lo_{\G_2}-\Lo_{\G_1}\right)\ \Ui_{\G_1}
= \Us_{\G_2}\,\Lo_{\Po_{H} (F_2 - F_1)}\,\Ui_{\G_1}\\
= \Us_{21}\,\Lo_{\,\Us_{\G_1}\Po_{H} (F_2 - F_1)}\,.
\]
Here, we have used the invariance of Poisson brackets by canonical transforms $\Lo_F\Ui_{\G}=\Ui_{\G}\Lo_{\Us_{\G}F}$. Because of the intertwining  relation \eref{Intereq},  this generator  is  always secular:
\[\G_{21}=\Us_{\G_1}\, \Po_{H} ( F_2 -  F_1) =
 \Po_{H_0} \Us_{\So_H {\partial H}/{\partial\,\ale} + \Po_H F_1}  \, ( F_2 -  F_1).
\]
Therefore, the secular normalised Hamiltonians obtained using different uniqueness conditions are related by the Lie-Deprit transform
with the  secular generator.
This corresponds
to  the Bruno theorem for generating functions  \citep{Bruno}. 

For  the non-resonance  system with incommensurable frequencies, all the integrals  commute.
As a consequence,
a non-resonance normalised Hamiltonian is unique \citep{Koseleff}.

This  is not so in the case of resonance. Because  resonance relations lead to non-commuting  integrals,  the 
normalised resonance Hamiltonian  depends on the uniqueness condition. We can obtain an explicitly secular expression 
 for the Hamiltonian following \citet{Vittot0}. 
 For the generator $\G$ given by \eref{Final}, consider the derivative:
\begin{align*}
 \frac{\partial}{\partial\,\ale}\tilde H &=  \left(\frac{\partial}{\partial\,\ale}  \Us_{\G}\right) H + 
 \Us_{\G} \frac{\partial\, }{\partial\,\ale} H =
\Us_{\G} \, \left(\Lo_\G H+\frac{\partial H }{\partial\,\ale} \right) \\
&=\Us_{\G} \, \left(\frac{\partial H }{\partial\,\ale}-\Lo_H \left(  \So_H \frac{\partial H}{\partial\,\ale} + \Po_H F \right) \right) 
=\Us_{\G} \, \left(   \Po_H \frac{\partial H}{\partial\,\ale} - \D_H F  \right) \\
&= \Us_{\G}  \, \Po_H \, \left(\frac{\partial H }{\partial\,\ale} - \D_H F\right) =
 \Po_{H_0}\, \Us_{\G} \, \left(\frac{\partial H }{\partial\,\ale}  - \D_H F\right).
\end{align*} 
Here, we  used \eref{depriteq}, the properties of the perturbed operators \eref{pertprop} and the intertwining  of  projectors \eref{Intereq}. Therefore,
\begin{equation}
\tilde H =  H_0 + 
\Po_{H_0} \int_0^\ale \Us_{\So_H {\partial H}/{\partial\,\ale} + \Po_H F}  \left(\frac{\partial H }{\partial\,\ale}  - \D_H F\right) \,\rmd\ale.\label{FinalH}
\end{equation}
This  secular expression demonstrates  the explicit dependence on  $\Po_H F$. However, it requires more computational resources than
the standard normalisation. 

\subsection{Gustavson integrals}

It is  interesting to find physically
meaningful quantities that  do not depend on an artificial choice of the uniqueness condition $\Po_H F({\bf x},\ale)$.
Our previous work \citep{nikolaev4} demonstrated that for linearly perturbed systems, the Gustavson integrals have such a property.
Now, we are able to  apply similar considerations to a general  case.

Consider   a system with constant unperturbed frequencies  that obey
${\textbf r}$ linearly independent resonance relations 
$(\vecp \omega,\vecp D_k)=0$, $k=1,\ldots ,{\textbf r}$. 
Such resonance relations result in additional non-commuting integrals. 
In the Birkhoff and action-angle representations, the centre of the corresponding  algebra 
consist of the ${\textbf d}-{\textbf r}$ quantities:
\[\tilde I_m = \sum_{j=1}^{\textbf d} \beta_{mj}\tilde\zeta_j \tilde\eta_j = (\vecp\beta_m, \vecp J),\qquad m=1,\ldots,{\textbf d}-{\textbf r},
\]
where ${\vecp \beta}_m$ is a set of ${\textbf d}-{\textbf r}$ independent vectors that are orthogonal to all the  
resonance vectors
$\vecp D_k$. 
These integrals  commute with all the integrals of the unperturbed system and, therefore 
with the normalised Hamiltonian ${\widetilde H}$ \citep{Gustavson}. 
In the operator notation, for any analytic function 	$\tilde F({\bf{x}})$:
\begin{align*}
&\tilde I_m = \Po_{H_0} \tilde I_m, \qquad\qquad m=1,\ldots,{\textbf d}-{\textbf r},\\
&[ \tilde I_m, \Po_{H_0} \tilde F ]=0,\qquad\\
&[ \tilde I_m,  {\widetilde H} ]=0.
\end{align*}
After the transform back to the initial variables,  the quantities
\[ I_m = \Ui_\G\, \tilde I_m,\qquad m=1,\ldots,{\textbf d}-{\textbf r},
\]
 become formal integrals  of the perturbed system. These Gustavson integrals belong to  the image of $\Po_{H}$
 and  commute with all the functions in this space: 
\begin{align*}
&I_m = \Ui_\G \, \Po_{H_0} \tilde I_m =  \Po_{H} \,  \Ui_\G \,\tilde I_m = \Po_{H} I_m, \qquad m=1,\ldots,{\textbf d}-{\textbf r},\\
&[ I_m, \Po_{H}  F ]= \Ui_\G \, [ \tilde I_m, \Po_{H_0}  \Us_\G \, F ] = 0,\\
&[  I_m,  H ]= \Ui_\G \, [ \tilde I_m, \tilde H ] = 0.
\end{align*}
Here, we have again used the intertwining  of the projectors  \eref{Intereq}.
Due to the above properties,  the derivative of $I_m({\bf{x}},\ale)$ does not depend on the $\Po_H F$ part of $\G$:
\begin{align*}
\frac{d}{d\,\ale}I_m (\ale)  &=  \left(\frac{\partial}{\partial\,\ale}   \Ui_\G\right) \tilde I_m  =
-\Lo_{\G}   \Ui_\G\, \tilde I_m = - \Lo_{\G} \,I_m =\\
&= - \Lo_{\So_H {\partial H}/{\partial\,\ale}} \, I_m - \Lo_{\Po_H F}\, I_m  = - \Lo_{\So_H {\partial H}/{\partial\,\ale}}\,  I_m .
\end{align*}
Therefore, the Gustavson formal integrals $I_i(\ale)$  are  insensitive to the  uniqueness condition.  
The corresponding series diverge \citep{Contopoulos2003}, but  are useful  for exploring 
the regions of regular dynamics.

The unperturbed Hamiltonian $H_0$  can be chosen as the initial function for the Gustavson integral.
The resulting series is  known as the Hori  formal first integral  \citep{Hori}:
\begin{align*}
I_H =&  \ale^{-1} \left( H - \Ui_{\So_H {\partial H}/{\partial\,\ale} } H_0 \right) = \Po H_1 + \ale \left( \Po H_2
- \So\Lo_{H_1}\Po H_1 - \frac{1}{2}\Po\Lo_{H_1}\So H_1 \right)\\
 &+ \ale^2 \left( \So \Lo_{H_1} \So \Lo_{H_1} \Po H_1 + \frac{1}{2} \So \Lo_{H_1} \Po \Lo_{H_1} \So H_1 \!
 +\frac{1}{3} \Po \Lo_{H_1} \So \Lo_{H_1} \So H_1\!
 - \!\frac{2}{3} \Po \Lo_{H_1} \So^2 \Lo_{H_1} \Po H_1 \right.\\
 &\left. - \frac{1}{3} \Po \Lo_{H_1} \Po \Lo_{H_1} \So^2 H_1  - \So \Lo_{H_1}\Po H_2 - \So\Lo_{H_2}\Po H_1 + \Po H_3 \right) 
+\Ord(\ale^{3}).
\end{align*}
It is also applicable to  systems with non-constant unperturbed frequencies.

\section{Computational aspects}
A major difference between this and  the classical canonical perturbation algorithms 
is the  explicit non-recursive formulae.  Usually, 
perturbation  computations normalise the Hamiltonian order by order.
In contrast, we    compute  the  generator  up to  the  
desired order directly. 
Next, the direct Lie-Deprit transform normalises the perturbed Hamiltonian,  and the inverse  transform computes its integrals.
\subsection{An explicit algorithm for the  generator}  
To normalise the Hamiltonian $H$ up to the order $\Ord{(\ale^{N+1})}$, we must compute the generator  $\G$ up to  $\Ord{(\ale^{N})}$.
It is possible, but not efficient,  to use the combinatorial sum  \eref{explicit}.
A more elegant algorithm  follows from the Neumann series \eref{Neumann}. 
Consider  the following expression:
\[\R_{H}(z) \frac{\partial H}{\partial \ale}=\R_{H_0}(z)\sum_{n=0}^{N-1}{\left(-\Lo_{V}\R_{H_0}(z)\right)}^n\ \frac{\partial H}{\partial \ale}
 +\Ord{(\ale^{N})}, \qquad\textrm{where}\ V = \sum_{k=1}^{N-1} \ale^k H_k.
\]
Its   coefficient of $z^0$   gives the generator $\G =  \So_{H}\,  \Ha$ up to the desired accuracy. 

In order to avoid  the negative powers of the ancillary variable $z$ during computation, we  use the operator:
\[ \Rz(z)=-\Po_{H_0} + \sum_{n=1}^N z^n \,\So_{H_0}^{n}\equiv z\,\R_{H_0}(z)+\Ord{(z^{N+1})}.
\]

In the Birkhoff representation, 
the action of this linear operator on a monomial is computed as follows: 
\footnotetext[1]{The rational expression ${z}/({\rmi\,(\vecp\omega,\vecp k - \vecp m) - z})\ \zeta^{\vecp k}\eta^{\vecp m}$ is preferable for   some computer algebra systems.}
\[
\Rz\,   \zeta^{\vecp k}\eta^{\vecp m} = 
\begin{cases} \qquad\qquad {\bf -}  \quad  \zeta^{\vecp k}\eta^{\vecp m},& \textrm{when}\  {(\vecp\omega,\vecp k - \vecp m) =0},\\
\sum\limits_{n=1}^N  \left(\frac{z}{\rmi\,(\vecp\omega,\vecp k - \vecp m)}\right)^n  \zeta^{\vecp k}\eta^{\vecp m}
, & 
\textrm{when} \  {(\vecp\omega,\vecp k - \vecp m) \ne 0}\footnotemark[1].
\end{cases}
\]
 The expression in the action-angle variables is similar. 

Our algorithm sequentially constructs  the   $N-1$ functions:
\[F_1({\bf x},\ale,z)=\Rz\,\frac{\partial H}{\partial \ale},\quad\ldots,\quad
F_{n}({\bf x},\ale,z)=-\Rz\,\Lo_{V}F_{n-1},\qquad n=1,\ldots,N.
\]
At each step, a computer algebra system  automatically discards all the terms 
containing $\ale^{N}$ and $z^{N+1}$.
The generator up to $\Ord{(\ale^{N})}$ is directly computed from these functions:
\[\G(\ale) =  \So_{H}\, \frac{\partial H}{\partial \ale} = \sum_{n=1}^{N} F_n[z^{n}],
\]
where $F_n[z^{n}]$ denotes the coefficient of $z^{n}$ in the function $F_n({\bf x},\ale,z)$.
\subsection{Computation of the Lie-Deprit transforms}
Traditionally, the Lie-Deprit transform $\Us$ and the generator $\G$ are  simultaneously computed order by order using the  triangular algorithm  \citep{Deprit}. 
In our approach, 
 we  compute the generator independently and can  use faster algorithms.
\paragraph{The inverse transform:} Let us first consider  the  Gustavson integrals $ I_m({\bf x},\ale)$. These quantities are computed by the application of $\Ui_\G$  to the quantities $\tilde I_m({\bf x})$ that do not depend on $\ale$. An efficient computation should utilize this property  \citep{Henrard72}.

The Henrard algorithm is based on the recursive relation \eref{Caryrel}.
In our particular case, it iteratively constructs  $N+1$ ancillary functions:
\[
\tilde {f}_{0}=\tilde I_m({\bf x}),
\quad\ldots,\quad
\tilde {f}_{n}=-\frac{1}{n}\sum_{k=0}^{n-1} \Lo_{\G_{n-k-1}} \tilde f_k,\qquad n=1,\ldots,N.
\]
These functions are the coefficients of the approximation for the Gustavson integral:
\[ I_m({\bf x},\ale)= \Ui_\G \tilde I_m =\sum_{n=0}^N \ale^n \tilde {f}_{n}+\Ord(\ale^{N+1}).\]

\paragraph{The direct transform:}
Inspired by the simplicity and the efficiency of the above computation, we propose a similar algorithm for the direct Lie-Deprit transform. The normalised Hamiltonian  
 can be written as the following sum:
\[
{\widetilde H}=  \Us_\G H = \left( \sum_{n=0}^N \ale^n\, \Us_n \right)  \left(\sum_{n=0}^N \ale^n  H_n\right) +\Ord(\ale^{N+1}) 
= \sum_{n=0}^N  \Us_n {f}{}_{n}^{(N)}+\Ord(\ale^{N+1}).
\]
Here, we  introduce the  ancillary double-indexed 
functions ${f}^{(N)}_{n}({\bf x},\ale) = \sum_{k=n}^{N} \ale^{k}  H_{k-n}$.

The   relations \eref{Caryrel} 
allows us repeatedly express   the operators $\Us_n$ by means of 
$\Us_{\text{less then }n}$:
\[
{\widetilde H}= 
 \Us_N{f}_{N}^{(N)}+ \sum_{k=0}^{N-1}  \Us_k {f}_{k}^{(N)} = 
 \sum_{k=0}^{N-1}  \Us_k {f}_{k}^{(N-1)}
 =\ldots= \sum_{k=0}^{n}  \Us_k {f}_{k}^{(n)}=\ldots ,
\]
where  the functions ${f}_{k}^{(n)}$, $n=N-1,\ldots,0$ are computed using the  relations:
\[
{f}_{k}^{(n)} =
 {f}_{k}^{(n+1)} + \frac{1}{n+1} \Lo_{\G_{n-k}} {f}_{n+1}^{(n+1)},
 \qquad  k=0,\ldots,n.
\]
 Finally, all the operators $\Us$ disappear and we obtain  the transformed Hamiltonian:
\[
{\widetilde H} =\ldots= {f}_{0}^{(1)} + \Us_1 {f}_{1}^{(1)}=
{f}_{0}^{(0)}.
\]
The  functions ${f}_{k}^{(n)}({\bf x},\ale)$ form a triangle similar to that of Deprit, but their computation 
does not include the summation. 

\section{Examples and comparison with other canonical algorithms}
In the following examples, we  compare   the  normalisation 
of the model Hamiltonians 
 using our explicit algorithm 
 with the results of the major canonical perturbation approaches,
 including \citet{Gustavson},  Hori-\citet{Mersman}, \citet{Henrard72}
and \citet{DragtFinn} algorithms.

\subsection{Pendulum}
The Hamiltonian of the standard one-dimensional pendulum   is $P^2/2+(1-\cos{Q})$.
The canonical non-univalent scale transform
$P=\sqrt{\ale}p$, $Q=\sqrt{\ale}q$ 
allows us to introduce the perturbation parameter:
\[ H=\frac{1}{2}( p^2+ q^2)+\sum_{k=1}^\infty \frac{(-1)^{k+1}}{(2(k+1))!} q^{2(k+1)}\ale^k .
\]
The  first orders of the normalising generator and
the normalised Hamiltonian are:
\begin{align*}
\G=& \frac{5}{192} p_1 q_1^3 + \frac{1}{64} p_1^3 q_1
+ \ale  \left( \frac{17}{7680} p_1 q_1^5 + \frac{1}{192} p_1^3 q_1^3 + \frac{1}{512} p_1^5 q_1 \right)
+\Ord(\ale^3),
\\
{\widetilde{H}}=& \frac{1}{2}\left(p^2+ q^2\right)
-\frac{1}{64} \left(p^2+ q^2\right)^2 \ale 
- \frac{1}{2048} \left(p^2+ q^2\right)^3 \ale^2
- \frac{5}{131072}\left(p^2+ q^2\right)^4 \ale^3\\
&{}- \frac{33}{8388608} \left(p^2+ q^2\right)^5 \ale^4
 - \frac{63}{134217728} \left(p^2+ q^2\right)^6 \ale^5
+\Ord(\ale^6). 
\end{align*}
This is the typical structure of the 
 perturbation series for non-resonance systems.
Although different canonical perturbation  methods build different near-identity  canonical normalising  transforms, 
the normalised  Hamiltonians of non-resonance   systems are identical  \citep{Koseleff}.

\subsection{Toda 2D system}
The Toda  two-dimensional system has the following Hamiltonian:
\[ H_{T}= 
\frac{1}{2} \left(P_1^2 + P_2^2\right)+ \frac{1}{24} \left( \rme^{2 Q_2+2 \sqrt{3} Q_1}
+\rme^{2 Q_2-2\sqrt{3} Q_1}+\rme^{-4 Q_2}\right)-\frac{1}{8}.
\]
Its integrability was first stated by 
 \citet{Ford73}  using perturbative method.
 
The scale transform
$P_j=\ale p_j$, $Q_j=\ale q_j$ 
introduces the perturbation parameter.
It is worth  mentioning that the  first two orders of an $\ale$ expansion for the Toda 2D system  coincide with the  Henon-Heiles Hamiltonian:
\begin{align*}
H=&\frac{1}{2} \left(p_1^2 + q_1^2 +p_2^2 + q_2^2\right)
+    \ale \left(  - \frac{1}{3}q_2^3 + q_1^2 q_2 \right)
       + \ale^2  \left(  \frac{1}{2} q_2^4 + q_1^2 q_2^2 +  \frac{1}{2} q_1^4 \right)\\
       &{}+ \ale^3  \left(  -  \frac{1}{3} q_2^5 +  \frac{2}{3} q_1^2 q_2^3 + q_1^4 q_2 \right)
       + \ale^4  \left(  \frac{11}{45} q_2^6 +  \frac{1}{3} q_1^2 q_2^4 + q_1^4 q_2^2 +  \frac{1}{5} q_1^6 \right)
       +\Ord(\ale^5).
\end{align*}
 Because this is a 1:1 resonance system, the corresponding  normalised  Hamiltonian contains  mixed terms, 
 such as $\zeta_1 \eta_2$:
\begin{align*}\widetilde H =&       \rmi\,(\zeta_2 \eta_2 +  \zeta_1 \eta_1)
       + \ale^2   \left( -\frac{1}{3} \zeta_2^2 \eta_2^2 + \frac{1}{3} \eta_1^2 \zeta_2^2 - \frac{4}{3}\zeta_1 \eta_1 \zeta_2 \eta_2 
 + \frac{1}{3}\zeta_1^2 \eta_2^2 - \frac{1}{3}\zeta_1^2 \eta_1^2 \right) \\
       &+ \rmi \ale^4 \left( -\frac{5}{27}  \zeta_2^3 \eta_2^3 - \frac{7}{9}  \eta_1^2 \zeta_2^3 \eta_2 
+ \zeta_1 \eta_1 \zeta_2^2 \eta_2^2 - \frac{7}{9} \zeta_1 \eta_1^3 \zeta_2^2 - \frac{7}{9} \zeta_1^2 \zeta_2 \eta_2^3 \right.\\
&\left.+  \zeta_1^2\eta_1^2\zeta_2\eta_2 
- \frac{7}{9} \zeta_1^3 \eta_1 \eta_2^2 - \frac{5}{27}\zeta_1^3 \eta_1^3 \right) + \Ord(\alpha^{6}).
\end{align*}
The following shows the  series for the Gustavson integral:
\begin{align*}I_H=& \ale^{-2} (H - \Ui_\G H_0 )  
=\frac{1}{12} q_2^4 +\frac{1}{6}p_2^2 q_2^2 + \frac{1}{12} p_2^4 + \frac{1}{6} q_1^2 q_2^2 + \frac{1}{2} q_1^2 p_2^2
          + \frac{1}{12} q_1^4 \\
&          {}- \frac{2}{3} p_1 q_1 p_2 q_2 + \frac{1}{2} p_1^2 q_2^2 + \frac{1}{6} p_1^2 p_2^2 
          + \frac{1}{6} p_1^2 q_1^2 + \frac{1}{12} p_1^4
+ \ale \left(  - \frac{1}{9}q_2^5 \right.\\
&\quad{}- \frac{1}{9} p_2^2 q_2^3 + \frac{2}{9} q_1^2 q_2^3 + \frac{5}{3} q_1^2 p_2^2 q_2
+ \frac{1}{3} q_1^4 q_2 - \frac{2}{3} p_1 q_1 p_2 q_2^2 + \frac{4}{3} p_1 q_1 p_2^3  \\
&\quad{}- \frac{2}{3} p_1 q_1^3 p_2- \frac{7}{9} p_1^2 q_2^3 - \frac{4}{3} p_1^2 p_2^2 q_2 + p_1^2 q_1^2 q_2 - \frac{4}{9} p_1^3 q_1 p_2
          \left.+ \frac{4}{9} p_1^4 q_2 \right)
+ \Ord(\alpha^{2}).
\end{align*}

Again,  the canonical perturbation  methods build  different near-identity  canonical $\ale$-parametrized normalising  transforms. 
The \citet{Gustavson} method uses the generating functions of mixed variables according to Poincar\'e-Birkhoff-Von Zeipel approach.
For the same purpose,  the Hori-\citet{Mersman} algorithm uses the  Lie series (Magnus expansion), while  \citet{DragtFinn} construct a product of operatorial exponents (Fer expansion).
The \citet{Deprit} method builds the direct Lie-Deprit transform for the Hamiltonian normalisation, but its uniqueness condition differs from our explicit formula \eref{exp1}.
Finally, the \citet{Henrard72} algorithm is based on  the inverse  Lie-Deprit transform. Traditionally, all these methods construct a non-secular generator.

All these normalising  transforms   
are equivalent to    the 
 Lie-Deprit transforms  with different uniqueness conditions $\Po_{H} F$.
Because of the non-commutativity of the unperturbed integrals, the  Hamiltonians of a resonance system 
normalised by different methods   differ from each other.

As can be seen from our supplemental demonstrations, 
the normalised Toda 2D Hamiltonians  
constructed by the Deprit, Hori, Gustavson, Dragt-Finn and Henrard methods using the non-secular uniqueness condition
  differ from each other starting from the  $8\textsuperscript{th}$ 
 order. 
 All of them are connected by secular canonical transforms. 
 
Moreover, we can either normalise the Hamiltonian in the $(p,q)$ variables or 
transform it first to the $(\eta,\zeta)$ variables, normalise it there and transform back.
The results of  the Lie algebraic algorithms are equal.
Since the generating functions are not invariant under  canonical transforms, the results of the Gustavson normalisations 
 in the $(p,q)$ and $(\eta,\zeta)$ variables 
differ starting from the  $6\textsuperscript{th}$  perturbation order.

As expected, the series for the Gustavson integrals for all the methods coincide up to the highest order that we computed. 

\paragraph{Henrard method:} It is interesting that the  Hamiltonian normalised by our explicit method  coincides with that  constructed by the \citet{Henrard72} method.
This is not occasional. The  Henrard method  uses the fact that  any direct Deprit transform $\Us_\G$ 
may be written as  an inverse
  Deprit transform $\Ui_{\widetilde{\G}}$ with the generator $\widetilde{\G}=-\Us_\G \G$. 
  
  The Henrard method normalises the perturbed Hamiltonian using this inverse Deprit transform because it is easier to compute. 
  It
   constructs the generator  $\widetilde{\G}(x,\ale)$ order by order using the requirement that $\widetilde H=\Ui_{\widetilde{\G}} H$ is secular. At each order it chooses the non-secular solution of the homological equation, therefore 
   $\Po_{H_0} \widetilde{\G} \equiv 0$. The  details can be found in  \citet{Koseleff}.

Due to the intertwining relations \eref{Intereq}:
\[\Po_{H_0} \widetilde{\G} = -\Po_{H_0} \Us_\G \G = - \Us_\G \Po_{H}  \G \equiv 0,
\]
 the corresponding Deprit generator obeys our  natural uniqueness condition $\Po_{H}  \G \equiv 0$.
  Thus, the  Henrard transform actually coincides with the Deprit transform with the generator obtained by the explicit formula
 \eref{exp1}.

\section{Computational efficiency}

Traditionally, the efficiencies of the perturbation methods are compared by the number of Poisson brackets \citep{Broer}.
However, this formal comparison does not take into account the resource consumption for solutions of homological equations,
series substitution and memory management. This is why we  
prefer a computer benchmarking of the particular realisations of algorithms. 

Unfortunately, even the benchmarking  efficiency is not unambiguous. 
High-order computations process large multi-gigabyte expressions containing millions of terms.
The performance of operations with such expressions depends
on  the computer algebra system used,    optimisations, the server CPU, RAM, OS,
 disks  and filesystem type.
Even the relative efficiencies of the methods may vary. 
Therefore, the following  results are only illustratory.
\begin{center}
\vbox{
\includegraphics[width=0.8\linewidth]{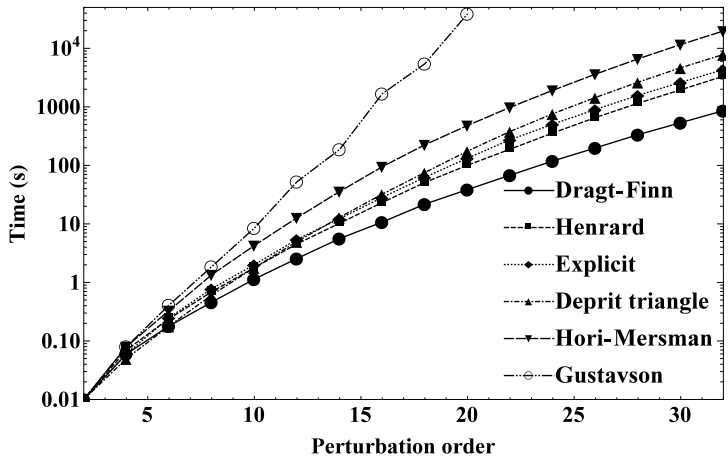}
\captionof {figure}{\label{fig:comp}H\'enon-Heiles system}
}
\end{center}

Figure \ref{fig:comp} compares 
the  times to compute the  normalised Hamiltonian and  the Gustavson integral for the Henon-Heiles system 
\[H_{HH}= \frac{1}{2} (p_1^2+q_1^2 + p_2^2+q_2^2)+ \alpha (q_1^2 q_2-\tfrac{1}{3} q_2^3)\]
on an Oracle\texttrademark{} Exadata X2-2 server with Intel Xeon X5675 (3.06 GHz) processor using Form~4.1 computer algebra system  \citep{Form}. 
The  Gustavson, Dragt-Finn and Henrad algorithms were run in the $(p,q)$ variables.
The other methods have the better  performance in the $(\eta,\zeta)$ variables.

We see that for a simple Hamiltonian with a limited number of perturbation terms 
only the Dragt-Finn and Henrad methods are faster than the explicit algorithm.

However, this changes for large Hamiltonians. For example, the number of  terms up to the $32\textsuperscript{nd}$ order  in the Toda 2D Hamiltonian 
exceeds 36000 in the $(\eta,\zeta)$ variables.
Corresponding computational times  are presented in Figure \ref{fig:comp2}. 
We see that our non-recursive explicit method becomes slower  than 
the recursive Lie algebraic  algorithms for high-order  normalisation of such  Hamiltonians.

\begin{center}
\vbox{\includegraphics[width=0.8\linewidth]{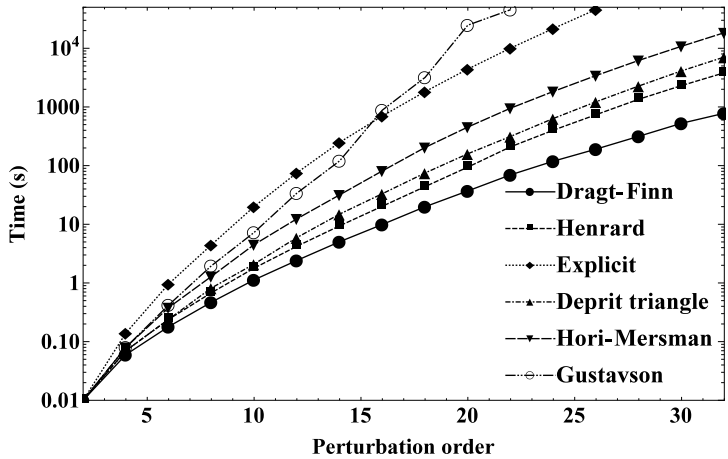}
\captionof {figure}{\label{fig:comp2}Toda 2D system}}
\end{center}

In all the cases we observed that the Dragt-Finn method had the superior performance. 
The Gustavson algorithm was the slowest  for high-orders computations.
Its   time  is dominated by the consecutive series substitutions during 
 the integral computation. 
 
It is worth  noting that the  computations were single-threaded. However, the  explicit algorithm
 can be readily parallelized and
made  scalable for contemporary multi-CPU and Cloud computing.

\section{Summary}
We have presented   the  application of the Kato  perturbation expansion  
   for the Laurent coefficients of the Liouvillian resolvent  to classical  Hamiltonians  
 represented by a power  series in a perturbation parameter.
This generalises our previous results \citep{nikolaev4} concerning 
the canonical intertwining of perturbed and unperturbed averaging operators and 
the  explicit expression  for the  generator of the normalising Lie-Deprit transform in any perturbation order. 

The approach allows for the systematic description of  non-uniqueness in the generator and normalised  Hamiltonian.
  We have also discussed the uniqueness of the Gustavson integrals and
  compared the  computational efficiency of the explicit expression 
  for  the H\'enon-Heiles and the Toda 2D systems
with the  efficiencies of classical canonical perturbation methods up to  the $32\textsuperscript{nd}$ perturbation order.

\begin{acknowledgements}
The author gratefully acknowledges  Professor S.~V.\ Klimenko,
Director of RDTEX Technical Support Center 
  S.~P.\ Misiura and V.~V.\ Romanova for encouragement and support.
\end{acknowledgements}

\bibliographystyle{spbasic}
\bibliography{kato_nlin}
\end{document}